\documentclass[12pt]{article}
\usepackage{latexsym,amsmath,amssymb}
\date{}
\usepackage{amsthm}
\allowdisplaybreaks

\newtheorem{theorem}{Theorem}

\newtheorem{lemma}[theorem]{Lemma}

\newtheorem*{conjecture*}{Conjecture}

\newcommand{\deq}{\stackrel{\scriptscriptstyle\triangle}{=}}

\def\D{\Delta}

\def\P{\mathbb P}
\def\E{\mathbb E}

\begin{document}
\title{A note on the random greedy triangle-packing algorithm}

\author{ {Tom Bohman\thanks{Department of Mathematical Sciences, Carnegie Mellon
University, Pittsburgh, PA 15213, USA. Email: {\tt tbohman@math.cmu.edu}.
Research supported in part by NSF grant DMS-0701183.}} \qquad
{Alan Frieze\thanks{Department of Mathematical Sciences, Carnegie Mellon
University, Pittsburgh, PA 15213, USA. Email: {\tt alan@random.math.cmu.edu}.
Research supported in part by NSF grant DMS-0721878.}} \qquad
{Eyal Lubetzky\thanks{Microsoft Research, One Microsoft Way, Redmond, WA 98052, USA. Email:
{\tt eyal@microsoft.com}.}}}


\maketitle

\begin{abstract}
The random greedy algorithm for constructing a large partial Steiner-Triple-System is defined as follows.
We begin with a complete graph on $n$ vertices and proceed to
remove the edges of triangles one at a time, where
each triangle removed is chosen uniformly at random
from the collection of all remaining triangles.
This stochastic process terminates once it arrives at a triangle-free graph.
In this note we show that with high probability the number of
edges in the final graph is at most $ O\big( n^{7/4}\log^{5/4}n \big) $.
\end{abstract}

\section{Introduction}

We consider the random greedy algorithm for triangle-packing.
This stochastic graph process begins with the
graph $ G(0) $, set to be the
complete graph on vertex set $[n]$, then proceeds to repeatedly remove the edges of randomly
chosen triangles (i.e.\ copies of
$K_3$) from the graph. Namely, letting $G(i)$ denote the graph that remains after $i$ triangles have been removed,
the $(i+1)$-th triangle removed is chosen uniformly at random from the set of all triangles in $ G(i)$.
The process terminates at a triangle-free graph
$G(M)$. In this work we study the random variable \(M\), i.e., the number of triangles removed until obtaining a triangle-free graph (or equivalently, how many edges there are in the final triangle-free graph).

This process and its variations play an important role in the
history of combinatorics.  Note that the collection of triangles removed
in the course of the process is a maximal collection of 3-element subsets of
\( [n] \) with the property that any pair of distinct triples in the collection have
pairwise intersection less than 2. 
For integers \( t < k < n  \) a {\em partial $(n,k,t)$-Steiner
system} is a collection of $k$-element subsets of an $n$-element set with the property that
any pairwise intersection of sets in the collection has cardinality less than $t$.
Note that
the number of sets in a partial $(n,k,t)$-Steiner system is at most \( \binom{n}{t}/ \binom{k}{t} \).
Let $ S(n,k,t) $ be the maximum number of $k$-sets in a partial $(n,k,t)$-Steiner system.
In the early 1960's Erd\H{o}s and Hanani~\cite{erdos} conjectured that for
any integers \(t < k\)
\begin{equation}
\label{eq:eh}
\lim_{n \to \infty } \frac{ S(n,k,t) \binom{k}{t} }{ \binom{n}{t} } = 1.
\end{equation}
In words, for any \( t < k \) there exist partial $(n,k,t)$-Steiner systems
that are essentially as large as allowed by the simple
volume upper bound.  This conjecture was proved by R\"odl~\cite{rodl} in the early 1980's by way of a randomized
construction that is now known as the R\"odl nibble.  This construction is a semi-random variation on the
random greedy triangle-packing process defined above, and thereafter such semi-random constructions have been successfully applied to establish various key results in Combinatorics over the last three decades (see e.g.~\cite{AS} for further details).

Despite the success of the R\"odl nibble, the limiting behavior of the random greedy
packing process remains unknown, even in the special case of triangle
packing 
considered here.  Recall that \(G(i)\) is the
graph remaining after \(i\) triangles have been removed.  Let \( E(i) \) be the edge
set of \( G(i) \).  Note that \( |E(i)| = \binom{n}{2} - 3i \) and that
\( E(M) \) is the number of edges in the triangle-free graph produced by the
process.  Observe that if we show \( |E(M)| = o(n^2) \) with non-vanishing probability then
we will establish~\eqref{eq:eh} for \( k=3, t=2\) and obtain that the
random greedy triangle-packing process produces an asymptotically optimal partial Steiner system.
This is in fact the case:  It was shown by Spencer~\cite{joel} and independently by R\"odl and Thoma~\cite{rodl} that \( |E(M)| = o(n^2) \) with high probability\footnote{Here and in what follows, ``with high probability'' (w.h.p.) denotes a probability tending to $1$ as $n\to\infty$.}. This was extended to \( |E(M)| \le n^{11/6 + o(1)} \)
by Grable in~\cite{grable}, where the author further sketched how similar arguments using more delicate calculations should extend to a bound of $n^{7/4+o(1)}$ w.h.p.

By comparison, it is widely believed that the graph produced by
the random greedy triangle-packing
process behaves similarly to the Erd\H{o}s-R\'enyi random graph with the same edge density, hence the
process should end once its number of remaining edges becomes comparable to the number of
triangles in the corresponding Erd\H{o}s-R\'enyi random graph.
\begin{conjecture*}[Folklore]
\label{conj:kahuna}
With high probability \( |E(M)|= n^{ 3/2 + o(1)} \).
\end{conjecture*}
\noindent
Joel Spencer has offered \$200 for a resolution of this question.

In this note we apply the differential-equation method 
to achieve an upper bound on \( E(M) \). In contrast to the 
aforementioned nibble-approach, whose application in this setting involves delicate calculations, 
our approach yields a short proof of the following best-known result:
\begin{theorem}
\label{thm:weak!} Consider the random greedy algorithm for triangle-packing on $n$ vertices. Let $M$ be the number of steps it takes the algorithm to terminate and let $E(M)$ be the edges of the resulting triangle-free graph. Then with high probability, $|E(M)| = O\big( n^{7/4}\log^{5/4} n \big)$.
\end{theorem}
\noindent
Wormald~\cite{nick2} also applied the differential-equation method to this problem, deriving 
an upper bound of $n^{2-\epsilon}$ on $E(M)$ for any $\epsilon < \epsilon_0 = 1/57$ while 
stating that ``some non-trivial modification would be required to equal or better Grable's result.''
Indeed, in a companion paper we combine the methods introduced here with some other ideas (and a 
significantly more involved analysis) to improve the exponent of the upper bound on $E(M)$ to 
about $1.65$. This follow-up work will appear in~\cite{bfl}.



\section{Evolution of the process in detail}

As is usual for applications of the differential equations method, we begin
by specifying the random variables that we track.  Of course, our main interest is in
the variable
\[ Q(i) \deq \text{ \# of triangles in } G(i)\,. \]
In order to track \( Q(i) \) we also consider the co-degrees in the graph $G(i)$:
\[ Y_{u,v}(i) \deq | \{ x \in [n] \,:\, xu, xv \in
E(i) \} | \]
for all \( \{u,v \} \in \binom{[n]}{2} \).  
Our interest in \( Y_{u,v} \) is motivated by the following observation:
If the $(i+1)$-th triangle taken is \( abc \) then
\[ Q(i+1) - Q(i) = Y_{a,b}(i) + Y_{b,c}(i) + Y_{a,c}(i) - 2\,. \]
Thus, bounds on \( Y_{u,v} \) yield
important information
about the underlying process.

Now that we have identified our variables, we determine the continuous
trajectories that they should follow.  We establish a correspondence with
continuous time by introducing a continuous variable $t$ and
setting \[ t = i / n^{2} \] 
(this is our time scaling).  We expect the graph
\( G(i) \) to resemble a uniformly chosen graph with \(n\) vertices and
\( \binom{n}{2} - 3i \) edges, which in turn resembles the Erd\H{o}s-R\'enyi graph \( G_{n,p} \)
with \[ p = 1 - 6i/n^2 = p(t)=1 - 6t\,. \]
(Note that we can view $p$ as either a continuous function of $t$ or
as a function of the discrete variable $i$.  We pass between
these interpretations of $p$ without comment.)
Following this intuition, we expect
to have \( Y_{u,v}(i) \approx p^2 n \) and \( Q(i) \approx p^3 n^3/ 6 \).  For ease of notation
define
\[ y(t) = p^2(t)~,\qquad  q(t) = p^3(t)/6\,.\]
We state our main result in terms of an error function that
slowly grows as the process evolves.  Define
\[ f(t) = 5 - 30 \log( 1-6t) = 5 - 30 \log p(t)\,.\]
Our main result is the following:
\begin{theorem}
\label{thm:starter}
With high probability we have
\begin{align}
&Q(i) \geq q(t) n^3 - \frac{ f^2(t) n^2\log n }{ p(t)} &\mbox{ and }  \label{eq:Qlower} \\
&\left| Y_{u,v}(i) - y(t) n \right| \leq  f(t) \sqrt{n \log n}  &\mbox{ for all $\{u,v\} \in \binom{[n]}{2}$}\,, \label{eq:Y}
\end{align}
holding for every
\[ i \le i_0=\tfrac16 n^2 -  \tfrac53 n^{7/4} \log^{5/4} n .\]
Furthermore, for all $i=1,\ldots,M$ we have
\begin{equation}
\label{eq:Qupper}
Q(i) \leq q(t) n^3 + \tfrac13 n^2 p(t)\,.
\end{equation}
\end{theorem}
\noindent Note that the error term in the upper bound \eqref{eq:Qupper}
{\em decreases} as the
process evolves.  This is not a common feature of applications of the differential equations method
for random graph process; indeed, the usual approach requires
an error bound that grows as the process evolves.  While novel techniques are introduced here
to get this `self-correcting' upper bound, two versions of `self-correcting'
estimates have appeared to date in applications of the differential equations
method in the literature (see \cite{mike} and \cite{nick}).  The stronger upper bound on the number of
edges in the graph produced by the random greedy triangle-packing process given in
the companion paper \cite{bfl} is proved by
establishing self-correcting estimates
for a large collection of variables (including the variable $Y_{u,v}$ introduced here).

Observe that~\eqref{eq:Qlower} (with $i=i_0$) establishes Theorem~\ref{thm:weak!}.
We conclude this section with a discussion of the implications of \eqref{eq:Qupper}
for the end of the process, the part of the process where there are fewer then \(n^{3/2} \) edges remaining.
Our first observation is that at any step \(i\) we can deduce a lower bound on the number
of edges in the final graph; in particular, for any \(i\) we have $ E(M) \ge E(i) - 3 Q(i)$.
We might hope to establish a lower bound on the number of edges remaining at the end of
the process by showing that there is a step \(i\)
where \( E(i) - 3Q(i) \) is large.  The bound~\eqref{eq:Qupper} is (just barely) too weak
for this argument to be useful.  But we can deduce the following.
Consider \( i = n^2/6 - \Theta( n^{3/2}) \); that is, consider
\(p = c n^{-1/2} \).  Once \(c\) is small enough the upper bound \eqref{eq:Qupper} is
dominated by the `error' term \( n^2p/3 \).  If $Q$ remains close to this upper bound then for the
rest of the process we are usually just choosing triangles in which every
edge is in exactly one triangle; in other words, the remaining graph is an approximate
partial Steiner triple system. If $Q$ drops significantly below this bound then
the process will soon terminate.

%

\section{Proof of Theorem~\ref{thm:starter}}
\label{sec:estimates}

The structure of the proof is as follows.
For each variable of interest and each bound (meaning both upper and lower) we introduce a
{\em critical interval}
that has one extreme at the bound we are trying to maintain and the other extreme slightly closer to the expected trajectory (relative to
the magnitude of the error bound in question).  The length of this interval is generally a function of $t$.  If a particular bound is violated
then sometime in the process the variable would have to `cross' this critical interval.  To show that this
event has low probability we introduce a collection of sequences of random variables, a sequence
starting at each step $j$ of the process.  This sequence stops as soon as the variable leaves the
critical interval (which in many cases would be immediately), and the sequence forms either a submartingale
or supermartinagle (depending on the type of bound in question).  The event that the bound in
question is violated is contained in the event that there is an index $j$ for which the corresponding
sub/super-martingale has a large deviation.  Each of
these large deviation events has very low probability, even in comparison with
the number of such events.  Theorem~\ref{thm:starter} then follows from the union bound.

For ease of notation we set
\[  i_0 = \tfrac16 n^2 - \tfrac53 n^{7/4}\log^{5/4}n ~,\qquad
p_0 = 10 n^{-1/4}\log^{5/4}n\,. \]
Let the stopping
time \(T\) be the minimum of \(M\) and the first step \( i <i_0\) at which \eqref{eq:Qlower} or \eqref{eq:Y} fail
and the first step \( i\) at which \eqref{eq:Qupper} fails.
Note that, since \( Y_{u,v} \) decreases
as the process evolves, if
\( i_0 \le i \le T \) then we have
\[ Y_{u,v} (i) = O \big( n^{1/2}\log^{5/2}n \big) \quad\mbox{for all $\{u,v\} \in \binom{[n]}{2}$}\,. \]

We begin with the bounds on \(Q(i) \).   The first observation is
that we can write the expected
one-step change in
\(Q\) as a function of $Q$.  To do this, we note that we have
\begin{equation}
\label{eq:starter}
\E[\D Q] = - \sum_{ xyz \in Q} \frac{ Y_{xy} + Y_{xz} + Y_{yz} - 2}{Q}
= 2 - \frac{1}{Q} \sum_{ xy \in E} Y_{xy}^2
\end{equation}
and
\[ 3Q = \sum_{xy \in E} Y_{xy}\,. \]
(And, of course, \( |E| = n^2p/2 - n/2 \).)  Observe that if $Q$ grows too large relative to
its expected trajectory then the expected change will be become more negative, introducing
a drift to $Q$ that brings it back toward the mean.  A similar phenomena occurs
if $Q$ gets too small.  Restricting our attention to a critical interval that is some distance
from the expected trajectory allows us to take full advantage of this effect.
This is the main idea in this analysis.

%
%

For the upper bound on $Q(i)$
our critical interval is
\begin{equation}
\label{eq:critQupper}
\left( q(t) n^3 + \tfrac14 n^2 p~,~ q(t) n^3 + \tfrac13 n^2p  \right).
\end{equation}
Suppose \( Q(i) \) falls in this interval.
%
Since Cauchy-Schwartz gives
\[ \sum_{xy \in E} Y_{xy}^2 \ge \frac{\left( \sum_{xy \in E} Y_{xy} \right)^2}{|E|} \ge \frac{9 Q^2}{ n^2 p/2}\,, \]
in this situation we have
\[ \E[ Q(i+1) - Q(i) \mid G(i) ] \le 2 - \frac{ 18 Q}{ n^2 p } < 2 - 3 n p^2 - \frac{9}{2} = -3 np^2 - \frac{5}{2}\,. \]
Now we consider a fixed index $j$.  (We are interested in those indices $j$ where $Q(j)$
has just entered the
critical window from below, but our analysis will formally apply to any $j$.)
We define the sequences of random variables \( X(j), X(j+1), \dots, X(T_j) \) where
\[ X(i) = Q(i) - q n^3 -  \frac{ n^2p}{3}  \]
and the stopping time \(T_j \) is the minimum of \( \max \{j,T\} \) and the smallest
index \( i \ge j \) such that \( Q(i) \) is not in the
critical interval \eqref{eq:critQupper}.  (Note that if $ Q(j) $ is not in the critical interval then we
have $ T_j = j $.)
In the event \( j \le i < T_j\) we have
\begin{equation*}
\begin{split}
\E[ X(i+1) - X(i) \mid G(i) ] & = \E[ Q(i+1) - Q(i) \mid G(i) ] - \left( q( t + 1/n^2) - q(t) \right) n^3 \\
& \hskip1cm - \left( p( t + 1/n^2) - p(t) \right) \frac{n^2}{3} \\
& \le - 3 n p^2 - \frac{5}{2} + 3 n p^2 + 2 + O \left(1/n \right) \\
& \le 0 \,.
\end{split}
\end{equation*}
So, our sequence of random variables is a supermartingale.  Note that if \( Q(i) \) crosses
the upper boundary in \eqref{eq:Qupper} at $i=T$ then, since the one step change in $Q(i)$ is
at most $3n$, there exists a step \( j \) such that
\[ X(j) \le - \frac{ n^2 p(t(j))}{4} + O(n) \] while \( T = T_j \) and \( X(T) \ge 0 \).
We apply Hoeffding-Azuma to bound the probability of such an event:
the number of steps is at most \( n^2 p( t(j))/6 \) and the
maximum 1-step difference is \( O( n^{1/2} \log^{5/2} n) \) (as \( i<T\) implies bounds on the co-degrees).
Thus the probability of such a large deviation beginning at step $j$ is at most
\[ \exp \left\{ - \Omega \left( \frac{  \left( n^2 p(t(j)) \right)^2 }{ (n^2 p( t(j))) \cdot  \big( n^{1/2} \log^{5/2}n   \big)^2 } \right)  \right\}
= \exp \left\{ - \Omega \left(  \frac{ n p( t(j))  }{ \log^5n} \right) \right\}\,.  \]
As there are at most $n^2$ possible values of $j$, we have the desired bound.

Now we turn to the lower bound on $Q$, namely \eqref{eq:Qlower}.  Here we work with the
critical interval
\begin{equation}
\label{eq:critQlower}
\left( q(t) n^3 - \frac{ f(t)^2 n^2\log n}{p}~,~  q(t) n^3 - \frac{ ( f(t) -1) f(t) n^2\log n }{p} \right)\,.
\end{equation}
Suppose $ Q(i) $ falls in this interval for some \( i < T \).
Note that our desired inequality is in the wrong direction for an application
of Cauchy Schwartz to \eqref{eq:starter}.  In its place
we use the control imposed on \( Y_{u,v}(i) \) by the condition $ i < T $.
For a fixed \( 3 Q = \sum_{uv \in E} Y_{u,v} \),
the sum \( \sum_{uv \in E} Y_{u,v}^2 \) is maximized when we make as many terms as large as possible.
Suppose this allows \( \alpha \) terms in the sum \( \sum_{ xy \in E} Y_{xy} \)
equal to \( np^2 + f \sqrt{n \log n} \)
and \( \alpha + \beta \) terms equal to \( n p^2 - f \sqrt{n \log n} \).  For ease of notation we
view \( \alpha, \beta \) as rationals, thereby allowing the terms in the maximum sum to split completely into these
two types.   Then we have
\[ \beta f  \sqrt{n \log n} =  |E| \cdot n p^2 -  3 Q = 3 q n^3 -3 Q - \frac{ n^2 p^2}{2} \,.  \]
Therefore, we have
\begin{equation*}
\begin{split}
\sum_{ xy \in E} &Y_{xy}^2  \le \alpha \left( n p^2 + f \sqrt{n \log n} \right)^2 +
( \alpha + \beta)\left( n p^2 - f\sqrt{n \log n} \right)^2 \\
& = \left( \frac{ n^2 p}{2} - \frac{n}{2} \right) \cdot n^2 p^4 +
\left( \frac{ n^2 p}{2} - \frac{n}{2} \right) \cdot f^2 n\log n - 2 \beta f p^2 n^{3/2}\log^{1/2} n  \\
& = \frac{ n^4 p^5}{2} + \frac{ f^2 p n^3 \log n}{2}
- \frac{1}{2}\left( n^3 p^4 + f^2 n^2 \log n\right)
- 2 p^2 n \left(3 q n^3 - 3 Q - \frac{n^2 p^2}{2} \right)  \\
& \le 6 n p^2 Q - \frac{ n^4 p^5}{2} + \frac{ f^2 p n^3\log n }{2} + \frac{ n^3 p^4}{2}\,.
\end{split}
\end{equation*}
Now, for $j < i_0$ define \( T_j\) to be the minimum of \(i_0\), \( \max \{j,T\} \) and the smallest index
\( i \ge j \) such that \( Q(i) \) is not in the critical interval \eqref{eq:critQlower}.
Set \[ X(i) = Q(i) - q(t)n^3 + \frac{ f(t)^2 n^{2}\log n }{ p(t)}\,. \]  For \( j \le i < T_j \) we have the bound
\begin{equation*}
\begin{split}
\E[ X(i+1) - X(i) &\mid G(i) ]   = \E[ Q(i+1) - Q(i) \mid G(i)] - n^3( q(t+1/n^2) - q(t)) \\
& \hskip1.5cm + \left( \frac{ f^2(t+1/n^2)}{ p( t + 1/n^2)} - \frac{ f^2(t)}{ p(t)} \right) n^2 \log n \\
& \ge 2 - 6 n p^2 + \frac{ n^4 p^5}{ 2 Q} - \frac{ f^2 p n^3 \log n}{2 Q } + O(p)  + 3 p^2 n + O(1/n)  \\
& \hskip1.5cm + \left( \frac{ 2 f^\prime f }{ p} + \frac{ 6 f^2}{ p^2} \right) \log n
+ O \left( \frac{ \log^3n}{ n^2 p^3} \right) \\
& \ge \frac{ ( f -1) f  n^2 \log n }{p} \cdot \frac{n^4 p^5}{ 2 (qn^3)^2} - \frac{ f^2 p n^3 \log n}{2 Q } + \left( \frac{ 2 f^\prime f }{ p} + \frac{ 6 f^2}{ p^2} \right)\log n  \\
& \ge \left[ \frac{ 18 f^2}{p^2} - \frac{ 18 f}{p^2} - (1+o(1)) \frac{ 3 f^2}{ p^2}  +  \left( \frac{ 2 f^\prime f }{ p} + \frac{ 6 f^2}{ p^2} \right) \right] \log n \\
&\ge 0\,.
\end{split}
\end{equation*}
If the process violates the bound \eqref{eq:Qlower} at step $T=i$ then there exists a $j < i$ such that
$ T = T_j$, $ X(T) = X(i) < 0 $ and
\[ X(j) >  \frac{ f( t(j)) n^2 \log n }{ p(t(j))} - O(n)\,. \]  The submartingale \( X(j), X(j+1), \dots X( T_j) \) has
length at most \( n^2 p(t(j))/6 \) and maximum one-step change \( O( n^{1/2} \log^{5/2} n ) \).  The probability that
we violate the lower bound \eqref{eq:Qlower} is at most
\[ n^2 \cdot \exp\left\{ - \Omega \left( \frac{ f^2(t(j)) n^4 \log^2 n/ p^2(t(j)) }{   n^2 p( t(j))\cdot n \log^5 n } \right) \right\} =
n^2 \cdot \exp\left\{ - \Omega \left( \frac{f^2( t(j)) n}{ \log^3 n } \right) \right\} = o(1)\,. \]

Finally, we turn to the co-degree estimate \( Y_{u,v} \).  Let $u,v$ be fixed.
We begin with the upper bound.  Our critical interval here is
\begin{equation}
\label{eq:critY}
\left( y(t) n + ( f(t)-5) \sqrt{n \log n}~,~ y(t) n + f(t) \sqrt{n \log n} \right)\,.
\end{equation}
For a fixed $j < i_0$ we consider the sequence of random
variables  $ Z_{u,v}(j), Z_{u,v}(j+1), \dots, Z_{u,v}(T_j) $ where
\[ Z_{u,v}(i) = Y_{u,v}(i) - y(t)n - f(t) \sqrt{n \log n} \]
and  \( T_j\) is defined to be the minimum of \(i_0\), \( \max \{j,T\} \) and the smallest index
\( i \ge j \) such that \( Y_{u,v}(i) \) is not in the critical interval \eqref{eq:critY}.
To see that this sequence forms a supermartingale, we
note that $i <T$ gives
\begin{equation*}
  \left | Q(i) -  q(t)n^3 \right| \leq \frac{f(t)^2  n^2\log n}{p(t)}\,,
\end{equation*}
%
%
%
and therefore
\begin{equation*}
\begin{split}
\E[ Z_{u,v}(i+1) - Z_{u,v}(i)] & \le -\sum_{x \in N(u) \cap N(v)} \frac{ Y_{u,x} + Y_{v,x} - 1_{uv \in E(i)}}{Q} \\
& \hskip1cm -n \left( y(t+1/n^2) - y(t) \right) - \sqrt{n \log n} \left( f(t+1/n^2) - f(t) \right) \\
& \le - \frac{  2 (y n + (f -5) \sqrt{n \log n})( y n - f
\sqrt{n \log n})}{ Q } + O \left( \frac{1}{ n^2p } \right) \\
& \hskip1.5cm  -  \frac{y^\prime(t)}{n} -  f^\prime(t) \frac{ \log^{1/2} n }{n^{3/2}}
+ O \left( \frac{1}{n^3 p^2} \right)  \\
& \le - \frac{  2 (y n + (f -5) \sqrt{n \log n} )( y n - f
\sqrt{n \log n})}{ q n^3 }  \\
& \hskip1cm   +  2 \cdot \frac{ f^2 n^2\log n }{p}  \cdot \frac{ (yn)^2}{ (q n^3)^2 }
-  \frac{y^\prime(t)}{n} -  f^\prime(t) \frac{ \log^{1/2}n }{n^{3/2}} + O \left( \frac{1}{ n^2p } \right) \\
& \le \frac{ 10 y  n^{3/2} \log^{1/2} n}{ q n^3} + \frac{ 14 f^2  n \log n}{ q n^3}  + O\left( \frac{ 1 }{ n^2 p} \right) -
f^\prime(t) \frac{ \log^{1/2} n }{n^{3/2}}
\end{split}
\end{equation*}
To get the supermartingale condition we consider each positive term here separately.
The following bounds would suffice
\[  \frac{ 60 }{ p } \le \frac{ f^\prime }{ 3 } ~,\qquad
\frac{ 84 f^2 \sqrt{\log n} }{ p^3 n^{1/2} } \le \frac{ f^\prime}{3 } ~,\qquad
\frac{ 1 }{ n ^{1/2} p } = o\left( f^\prime \sqrt{\log n} \right) \,. \]
The first term requires
\[f^\prime(t) \ge  \frac{ 180}{p(t)} =  \frac{180}{ 1 - 6t}\,. \]
We see that this requirement, together with the initial condition $ f(0) \ge  5 $, imposes
\[ f(t) \ge  5 -30\log( 1 - 6t) = 5 - 30 \log p(t)\,. \]
But this value for $f$ also suffices to handle the remaining terms as we restrict our
attention to \( p  \ge  p_0 = 10 n^{-1/4}\log^{5/4}n  \).
%
%
Thus, we have established that \( Z_{u,v}(i) \) is a supermartingale.

To bound the probability of a large deviation we recall a Lemma from \cite{r3t}.
A sequence of random variables $ X_0, X_1, \dots $ is {\em $ (\eta,N) $-bounded} if for all $i$ we have
\[ -\eta < X_{i+1}- X_i < N \,.\]

\begin{lemma}
Suppose \( 0 \equiv X_0, X_1, \dots \) is an
$( \eta, N )$-bounded submartingale for some $ \eta < N/10 $. Then for any $ a < \eta m $ we have $ \P( X_m < -a ) < \exp\big( - a^2 / ( 3 \eta N m ) \big)$.
\end{lemma}
\noindent
As \( -Z_{u,v}(j), -Z_{u,v}(j+1), \dots \) is
a $ (6/n , 2)$-bounded submartingale, the probability that we have
$ T = T_j$ with $ Y_{u,v}(T) > y n + f \sqrt{n \log n} $ is at most
\[ \exp \left\{ - \frac{ 25 n  \log n }{ 3 \cdot (6/n) \cdot 2 \cdot ( p(t(j)) n^2/6 )} \right\}
= \exp \left\{ - \frac{ 25 \log n}{6} \right\}
\,. \]
Note that there are at most $n^4$ choices for $j$ and the pair $u,v$.  As the argument for the
lower bound in \eqref{eq:Y} is the symmetric analogue of the reasoning we have
just completed, Theorem~\ref{thm:starter} follows. \qed


\end{document}